\documentstyle[amscd,amssymb,verbatim,12pt]{amsart}

\theoremstyle{plain}
\newtheorem{Prop}{Proposition}[section]
\newtheorem{Thm}[Prop]{Theorem}
\newtheorem{Cor}[Prop]{Corollary}

\newtheorem{Lem}[Prop]{Lemma}

\theoremstyle{definition}
\newtheorem{Def}[Prop]{Definition}

\theoremstyle{remark}
\newtheorem{Rem}[Prop]{Remark}

\def\int{\mathop{\roman{int}}}

\def\1{^{-1}}

\def\St{\text{St}}

\def\BB{{\mathcal B}}
\def\CC{{\mathcal C}}

\def\ZZ{{\mathbb Z}}
\def\dokaz{{\bf Proof. }}
\numberwithin{equation}{section}

\begin{document}
\title[
Coarse structures and group actions]%
   {Coarse structures and group actions
}

\author{N.~Brodskiy}
\address{University of Tennessee, Knoxville, TN 37996, USA}
\email{brodskiy@@math.utk.edu}

\author{J.~Dydak}
\address{University of Tennessee, Knoxville, TN 37996, USA}
\email{dydak@@math.utk.edu}
\author{A.~Mitra}
\address{University of Tennessee, Knoxville, TN 37996, USA}
\email{ajmitra@@math.utk.edu}

\date{ July 08, 2006}

\thanks{ The second-named author was partially supported
by Grant No.2004047  from the United States-Israel Binational Science
Foundation (BSF),  Jerusalem, Israel.
}

\begin{abstract}
The main results of the paper are:
\begin{Prop}\label{GenSvarc-Milnor}
A group $G$ acting coarsely on a coarse space $(X,\CC)$ induces a
coarse equivalence $g\to g\cdot x_0$ from $G$ to $X$ for any
$x_0\in X$.
\end{Prop}

\begin{Thm}\label{GenGromovThm}
Two coarse structures $\CC_1$ and $\CC_2$ on the same set
$X$ are equivalent if the following conditions are satisfied:
\begin{enumerate}
\item Bounded sets in $\CC_1$ are identical with bounded sets in $\CC_2$,
\item There is a coarse action $\phi_1$ of a group $G_1$
on $(X,\CC_1)$ and a coarse action $\phi_2$ of a group $G_2$ on $(X,\CC_2)$
such that $\phi_1$ commutes with $\phi_2$.
\end{enumerate}
\end{Thm}

They generalize the following two basic results of coarse
geometry:

\begin{Prop}[\v Svarc-Milnor Lemma~{\cite[Theorem 1.18]{Roe lectures}}] \label{Svarc-Milnor}
A group $G$ acting properly and cocompactly via isometries on a
length space $X$ is finitely generated and induces a
quasi-isometry equivalence $g\to g\cdot x_0$ from $G$ to $X$ for
any $x_0\in X$.
\end{Prop}

\begin{Thm}[Gromov~{\cite[page 6]{Gro asym invar}}] \label{GromovThm}
Two finitely generated groups $G$ and $H$ are quasi-isometric if
and only if there is a locally compact space $X$ admitting proper
and cocompact actions of both $G$ and $H$ that commute.
\end{Thm}

\end{abstract}

\maketitle


\section{Introduction}\label{section Introduction}

The proof in \cite{BDM} of the \v Svarc-Milnor Lemma was based on
the idea that isometric actions of groups ought to induce a coarse
structure on the group under reasonable conditions. Since left
coarse structures on countable groups are unique (in the sense of
independence on the left-invariant proper metric), \v Svarc-Milnor
Lemma follows.

In this paper we investigate cases where group actions on sets
induce a natural coarse structure on the set. As usual, the
uniqueness of the coarse structure is of interest.

We will use two approaches to coarse structures on a set $X$:
\begin{enumerate}
\item The original one of Roe \cite{Roe lectures} based on controlled
subsets of $X\times X$.
\item The one from \cite{DydHof} based on uniformly bounded families in $X$.
\end{enumerate}

The reason is that certain concepts and results have a more natural
meaning in a particular approach to coarse structures.
Recall that one can switch from one approach to another using
the following basic facts (see \cite{DydHof}):
\begin{itemize}
\item[a.] If $\{B_s\}_{s\in S}$ is uniformly bounded,
then $\bigcup\limits_{s\in S}B_s\times B_s$ is controlled.
\item[b.] If $E$ is controlled, then there is a uniformly bounded family
$\{B_s\}_{s\in S}$ such that $E\subset \bigcup\limits_{s\in S}B_s\times B_s$.
\end{itemize}

To define a coarse structure using uniformly bounded families
one needs to verify the following conditions
\begin{enumerate}
\item $\BB_1$ is uniformly bounded implies $\BB_2$
is uniformly bounded if each element of $\BB_2$
consisting of more than one point
is contained in some element of $\BB_1$.
\item $\BB_1,\BB_2$ uniformly bounded implies $\St(\BB_1,\BB_2)$
is uniformly bounded.
\end{enumerate}

\begin{Def}\label{BorDef}
A function $f\colon (X,\CC_X)\to (Y,\CC_Y)$ of coarse spaces is
{\it large scale uniform} (or {\it bornologous}) if $f(\BB)\in
\CC_Y$ for every $\BB\in\CC_X$.

$f$ is {\it coarsely proper} if $f^{-1}(U)$ is bounded for every
bounded subset $U$ of $Y$.

$f$ is {\it coarse} if it is large scale uniform and coarsely
proper.
\end{Def}

Recall that two functions $f,g\colon S\to (X,\CC_X)$ from a set
$S$ to a coarse space $(X,\CC_X)$ are {\it close} if the family
$\{\{f(s),g(s)\}\}_{s\in S}$ is bounded.

\begin{Def}
A coarse function $f\colon (X,\CC_X)\to (Y,\CC_Y)$ of coarse
spaces is a {\it coarse equivalence} if there is a coarse function
$g\colon (Y,\CC_Y)\to (X,\CC_X)$ such that $f\circ g$ is close to
$id_Y$ and $g\circ f$ is close to $id_X$.
\end{Def}

Here is a simple criterion for being a coarse equivalence
using the approach of \cite{DydHof}:

\begin{Lem} \label{CharOfCoarseEq}
A surjective coarse function $f\colon (X,\CC_X)\to (Y,\CC_Y)$ of
coarse spaces is a coarse equivalence if and only if $f^{-1}(\BB)$
is a uniformly bounded family in $X$ for each uniformly bounded
family $\BB$ in $Y$.
\end{Lem}

\dokaz Let $g\colon Y\to X$ be a selection for $y\to f^{-1}(y)$.
Put $\BB'=\{f^{-1}(y)\}_{y\in Y}\in\CC_X$.

If $g\colon (Y,\CC_Y)\to (X,\CC_X)$ is bornologous, then
$f^{-1}(\BB)$ refines $\St(g(\BB),\BB')$, resulting in
$f^{-1}(\BB)$ being uniformly bounded.

Let us show $g$ is bornologous if $f$ is a coarse equivalence.
Choose $h\colon (Y,\CC_Y)\to (X,\CC_X)$ that is bornologous and
$h\circ f$ is $\BB_1$-close to $id_X$ for some $\BB_1\in\CC_X$.
Therefore $h=h\circ f\circ g$ is $\BB_1$-close to $g$ and $g$ is
bornologous.

Assume $f^{-1}(\BB)$ is a uniformly bounded family in $X$ for each
uniformly bounded family $\BB$ in $Y$. If $g$ is bornologous, then
$f$ is a coarse equivalence as $f\circ g=id_Y$ and $g\circ f$ is
$\BB'$-close to $id_X$. If $\BB\in\CC_Y$, then $g(\BB)$ refines
$f^{-1}(\BB)$, so it is uniformly bounded and $g$ is bornologous.
\hfill $\blacksquare$

\begin{Cor}\label{EqualityOfCS}
Suppose $f\colon X\to Y$ is a surjective function
and $\CC_1$, $\CC_2$ are two coarse structures on $Y$.
If $\CC_X$ is a coarse structure on $X$
such that both $f\colon (X,\CC_X)\to (Y,\CC_i)$, $i=1,2$,
are coarse equivalences, then $\CC_1=\CC_2$.
 \end{Cor}

\dokaz Suppose $\BB\in\CC_1$ is uniformly bounded. Since
$f^{-1}(\BB)\in\CC_X$ by \ref{CharOfCoarseEq} and $f\colon
(X,\CC_X)\to (Y,\CC_2)$ is bornologous,
$\BB=f(f^{-1}(\BB))\in\CC_2$. Similarly, $\CC_2\subset \CC_1$.
\hfill $\blacksquare$

\begin{Rem}
We will see in~\ref{MonoFromZtoDih} that $f$ being surjective is
necessary.
 \end{Rem}

\section{Coarse structures on groups}\label{section CSGroups}

Given a group $G$ one can equip it with either the left coarse
structure $\CC_l(G)$ or right coarse structure $\CC_r(G)$. For
countable groups $G$ those structures are metrizable by proper
left-invariant (proper right-invariant) metrics on $G$.

In terms of controlled sets $E\in\CC_l(G)$ if and only if there is
a finite subset $F$ of $G$ such that $x^{-1}\cdot y\in F$ for all
$(x,y)\in E$. Similarly, $E\in\CC_r(G)$ if and only if there is a
finite subset $F$ of $G$ such that $x\cdot y^{-1}\in F$ for all
$(x,y)\in E$. Notice all functions $x\to g\cdot x$ ($g\in G$ being
fixed) are coarse self-equivalences of $(G,\CC_l(G))$ and all
functions $x\to x\cdot g$ are coarse self-equivalences of
$(G,\CC_r(G))$. We will primarily deal with the structure
$\CC_l(G)$ (notice $x\to x^{-1}$ induces isomorphism of structures
$\CC_l(G)$ and $\CC_r(G)$) but first we will characterize cases
where the two structures are identical.

\begin{Prop}\label{ComparingC_lAndC_r}
The following conditions are equivalent for any group $G$:
\begin{enumerate}
\item $\CC_l(G)=\CC_r(G)$,
\item $\CC_l(G)\subset \CC_r(G)$,
\item $\CC_r(G)\subset \CC_l(G)$,
\item $G$ is an FC-group (conjugacy classes of all elements are finite).
\end{enumerate}
 \end{Prop}

\dokaz $(3)\implies(4)$. Fix $a\in G$ and consider the family
$\{\{x,a\cdot x\}\}_{x\in G}$. It is uniformly bounded in
$\CC_r(G)$, so it must be uniformly bounded in $\CC_l(G)$ but that
means the set $\{x^{-1}\cdot a\cdot x\}_{x\in G}$ is finite, i.e.
the set of conjugacy classes of $a$ is finite. The same proof
shows $(2)\implies(4)$.

$(4)\implies(1)$. Given a uniformly bounded family $\BB$ in
$\CC_l(G)$ there is a finite subset $F$ of $G$ such that
$u^{-1}\cdot v\in F$ for all $u,v$ belonging to the same element
of $\BB$. Let $E$ be the set of conjugacy classes of all elements
of $F$. If $u,v$ belong to the same element of $\BB$, then there
is $f\in F$ so that $u^{-1}\cdot v=f$. Thus $v=u\cdot f$ and
$v\cdot u^{-1}=u\cdot f\cdot u^{-1}\in E$. Thus $\BB$ is uniformly
bounded in $\CC_r(G)$. The same argument shows $\CC_r(G)\subset
\CC_l(G)$. \hfill $\blacksquare$

\begin{Cor}\label{MonoFromZtoDih}
There is a monomorphism $i\colon\ZZ\to Dih_\infty$ from integers
to the infinite dihedral group $Dih_\infty$ that induces coarse
equivalences for both left coarse structures and the right coarse
structures but $\CC_l(Dih_\infty)\ne\CC_r(Dih_\infty)$.
\end{Cor}

\dokaz Consider the presentation $\{x,t\mid t^{-1}xt=x^{-1} \text{
and } t^2=1\}$ of $Dih_\infty$. Identify $\ZZ$ with the subgroup
of $Dih_\infty$ generated by $x$. Notice $\ZZ$ is of index $2$ in
$Dih_\infty$, so $\ZZ\to Dih_\infty$ is a coarse equivalence for
both left and right coarse structures. Since $\ZZ$ is Abelian,
those coincide on that group but
$\CC_l(Dih_\infty)\ne\CC_r(Dih_\infty)$ as the conjugacy class of
$x$ equals $\ZZ$. \hfill $\blacksquare$

\begin{Prop}\label{MultBornImpliesBothStructuresTheSame}
The multiplication $m:(G\times G,\CC_l(G)\times \CC_l(G))\to (G,\CC_l(G))$ is
large scale uniform if and only if  $\CC_l(G)=\CC_r(G)$.
\end{Prop}

\dokaz Suppose $F$ is a finite subset of $G$. Consider the
uniformly bounded family $\{F\times \{x\}\}_{x\in G}$ in
$\CC_l(G)\times \CC_l(G)$. Since $m(F\times \{x\})=F\cdot x$, the
family $\{F\cdot x\}_{x\in G}\in\CC(G)_l$. Thus $\CC_r(G)\subset
\CC_l(G)$ and $\CC_l(G)=\CC_r(G)$ by \ref{ComparingC_lAndC_r}.

Suppose $\CC_l(G)=\CC_r(G)$. It suffices to show that $\{m(x\cdot
F\times y\cdot E)\}_{(x,y)\in G\times G}$ is uniformly bounded for
every finite subsets $F$ and $E$ of $G$. Choose a finite subset
$E'$ and a function $f\colon G\to G$ such that $x\cdot E\subset
E'\cdot f(x)$ for all $x\in G$. Pick a finite subset $F'$ of $G$
and a function $g\colon G\to G$ such that $F\cdot E'\cdot y\subset
g(y)\cdot F'$ for all $y\in G$. Now $m(x\cdot F\times y\cdot
E)\subset x\cdot F\cdot E'\cdot f(y)\subset x\cdot g(f(y))\cdot
F'$ and the proof is completed. \hfill $\blacksquare$

\section{Inducing coarse structures by group actions}\label{section Inducing}

Our first task is to discuss cases of group actions of a group $G$
on a set $X$ inducing coarse structure $\CC_G$ on $X$ such that
$g\to g\cdot x_0$ is a coarse equivalence from $(G,\CC_l(G))$ to
$(X,\CC_G)$ for all $x_0\in X$.

\begin{Prop}\label{CreateFirstCS}
Suppose a group
$G$ acts transitively on a set $X$.
\begin{enumerate}
\item If there is a coarse structure $\CC_G$ on $X$ so that $g\to
g\cdot x_0$ is a coarse equivalence from $(G,\CC_l(G))$ to
$(X,\CC_G)$, then the stabilizer of $x_0$ is finite. \item If the
stabilizer of $x_0$ is finite, then there is a unique coarse
structure $\CC_G$ on $X$ so that $g\to g\cdot x_0$ is large scale
uniform. In that case $g\to g\cdot x_0$ is a coarse equivalence
from $(G,\CC_l(G))$ to $(X,\CC_G)$.
\end{enumerate}
\end{Prop}

\dokaz (1).  If $\gamma\colon g\to g\cdot x_0$ is a coarse
equivalence, then $\gamma^{-1}(x_0)$ must be bounded in $G$, i.e.
finite. Notice that $\gamma^{-1}(x_0)$ is precisely the stabilizer
of $x_0$.

(2). Assume the stabilizer $S$ of $x_0$ is finite. Define $\CC_G$
as follows: $\BB\in\CC_G$ if $\gamma^{-1}(\BB)$ is uniformly
bounded in $\CC_l(G)$. If $\CC_G$ is a coarse structure and
$\gamma\colon (G,\CC_l(G))\to (X,\CC_G)$ is bornologous, then
\ref{CharOfCoarseEq} says $\gamma$ is a coarse equivalence and the
uniqueness of $\CC_G$ follows from \ref{EqualityOfCS}.

Since
$\gamma^{-1}(\St(\BB_1,\BB_2))=\St(\gamma^{-1}(\BB_1),\gamma^{-1}(\BB_2))$,
$\BB_1,\BB_2\in\CC_G$ implies $\St(\BB_1,\BB_2)\in\CC_G$. Given
$\BB\in\CC_G$ we need to check that any family $\BB'$, whose
elements containing more than one point refine $\BB$, also belongs
to $\CC_G$. There is a finite subset $F$ of $G$ such that
$\gamma^{-1}(\BB)$ refines the family $\{g\cdot F\}_{g\in G}$. Put
$E=F\cup S$. If $\{x\}\in\BB'$, then $\gamma^{-1}(x)=h\cdot S$,
where $h\in G$ satisfies $x=h\cdot x_0$. Thus $\gamma^{-1}(\BB')$
refines $\{g\cdot E\}_{g\in G}$, so $\BB'\in\CC_G$. \hfill
$\blacksquare$

\begin{Prop}\label{CreateSecondCS}
Suppose a group $G$ acts on a set $X$. If there is a subset $U$ of
$X$ such that $X=G\cdot U$ and the stabilizer
$$S_U=\{g\in G\mid U\cap (g\cdot U)\ne\emptyset\}$$
of $U$ is finite, then there is a coarse structure $\CC_G$ on $X$
so that $g\to g\cdot x_0$ is a coarse equivalence from
$(G,\CC_l(G))$ to $(X,\CC_G)$ for all $x_0\in X$.
 \end{Prop}

\dokaz First define the bounded sets of $\CC_G$. Those are subsets
of sets of the form $F\cdot U$, where $F$ is any finite subset of
$G$. Second, define $\CC_G$ as families $\BB$ such that there is a
bounded set $V$ so that $\BB$ refines $\{g\cdot V\}_{g\in G}$.
Notice that, if $\BB'$ is a family whose elements containing more
than one point refine $\BB$, then $\BB'$ refines $\{g\cdot (V\cup
U)\}_{g\in G}$ and $V\cup U$ is bounded. Thus $\BB'\in\CC_G$.

The important property of bounded sets $V$ is that their
stabilizers $S_V=\{g\in G\mid V\cap (g\cdot V)\ne\emptyset\}$ are
finite. It suffices to prove that for $V=F\cdot U$, $F\subset G$
being finite. If $V\cap (g\cdot V)\ne\emptyset$, then there exist
elements $f_i\in F$, $i=1,2$, such that $(f_1\cdot U)\cap (g\cdot
f_2\cdot U)\ne\emptyset$ which implies $U\cap (f_1^{-1}gf_2\cdot
U)\ne\emptyset$. Thus $f_1^{-1}gf_2\in S_U$ and $g\in F\cdot
S_U\cdot  F^{-1}$ which proves $S_V$ is finite.

The second useful observation is that  $\St(V,\BB)$ is bounded for
any bounded set $V$ and any $\BB\in\CC_G$. Indeed, if $\BB$
refines $\{g\cdot W\}_{g\in G}$ for some bounded $W$, we may
assume $V\subset W$ in which case $V$ intersects only finitely
many elements of $\{g\cdot W\}_{g\in G}$. Since those are all
bounded and a finite union of bounded sets is bounded, we are
done.

Suppose $\BB_1,\BB_2\in\CC_G$ and choose bounded sets $V_i$,
$i=1,2$, such that $\BB_i$ refines $\{g\cdot V_i\}_{g\in G}$. Put
$V= \St(V_1,\{g\cdot V_2\}_{g\in G})$ and notice $V$ is bounded.
Our aim is to show $\St(\BB_1,\BB_2)$ refines $\{g\cdot V\}_{g\in
G}$. If $(h\cdot V_1)\cap (g\cdot V_2)\ne\emptyset$, then $V_1\cap
(h^{-1}\cdot g\cdot V_2)\ne\emptyset$, so $V_1\cup  (h^{-1}\cdot
g\cdot V_2)\subset V$ resulting in $\St(h\cdot V_1,\BB_2)\subset
h\cdot V$. \hfill $\blacksquare$

Let us point out that, surprisingly, the structure $\CC_G$ in
\ref{CreateSecondCS} does not have to be unique contrary to
typical categorical intuition.

\begin{Prop}\label{TwoCSExistOnDih}
There is an action of integers $\ZZ$ on the infinite dihedral
group $Dih_\infty$ such that $g\to g\cdot x_0$ are coarse
equivalences for both left and right coarse structures but
$\CC_l(Dih_\infty)\ne\CC_r(Dih_\infty)$.
 \end{Prop}

\dokaz Consider the presentation $\{x,t\mid t^{-1}xt=x^{-1} \text{
and } t^2=1\}$ of $Dih_\infty$. Identify $\ZZ$ with the subgroup
of $Dih_\infty$ generated by $x$. Notice $\ZZ$ is of index $2$ in
$Dih_\infty$, so $\ZZ\to Dih_\infty$ is a coarse equivalence for
both left and right coarse structures. Since $\ZZ$ is Abelian,
those coincide on that group but \ref{MonoFromZtoDih} says that
$\CC_l(Dih_\infty)\ne\CC_r(Dih_\infty)$. \hfill $\blacksquare$

\section{Actions by uniformly bornologous functions}\label{section ActionsbyULSU}

We want to generalize isometric actions to the framework of coarse
geometry. The appropriate concept is not only to require that each
function $x\to g\cdot x$ is bornologous but that those functions
are uniformly bornologous.

\begin{Def}\label{ActionByUBor}
A group $G$ acts on a coarse space $(X,\CC_X)$ by {\it uniformly
bornologous functions} if for any controlled set $E$ there is a
controlled set $E'$ such that $(g\cdot x,g\cdot y)\in E'$ for all
$(x,y)\in E$ and all $g\in G$.
\end{Def}

\begin{Prop}
A group $G$ acts on a coarse space $X$ by {\it uniformly
bornologous functions} if and only if for any uniformly bounded
family $\BB=\{B_s\}_{s\in S}$ in $X$ the family
$G\cdot\BB=\{g\cdot B_s\}_{(g,s)\in G\times S}$ is uniformly
bounded.
\end{Prop}

\dokaz Suppose the action is by uniformly bornologous functions
and $\BB=\{B_s\}_{s\in S}$ is a uniformly bounded family. Put
$E=\bigcup\limits_{s\in S}B_s\times B_s$ and notice it is a
controlled set. Pick a controlled set $E'$ such that $(g\cdot
x,g\cdot y)\in E'$ for all $g\in G$ and all $(x,y)\in E$. Define
$\BB'$ as the family of all $B\subset X$ satisfying $B\times
B\subset E'$. It is a uniformly bounded family containing $G\cdot
\BB$.

Suppose the family $G\cdot\BB=\{g\cdot B_s\}_{(g,s)\in G\times S}$
is uniformly bounded for any uniformly bounded family
$\BB=\{B_s\}_{s\in S}$ in $X$. Assume $E$ is a symmetric
controlled set containing the diagonal. Consider the family $\BB$
of all sets $B\subset X$ such that $B\times B\subset E\circ E\circ
E\circ E$ and let $E'=\bigcup\limits_{B\in \BB, g\in G}g\cdot
B\times g\cdot B$. It is a controlled set and, if $(x,y)\in E$,
then $\{x,y\}\times \{x,y\}\subset E\circ E\circ E\circ E$ and
$(g\cdot x,g\cdot y)\in E'$. \hfill $\blacksquare$

\begin{Cor}\label{BornImpliesUBorn}
Let $G$ be a group and $(X,\CC_X)$ be a coarse space. If
$\phi:(G\times X,\CC_l(G)\times \CC_X)\to (X,\CC_X)$ is
bornologous, then the action of $G$ on $(X,\CC_X)$ is by uniformly
bornologous functions.
\end{Cor}

\dokaz Given a uniformly bounded family $\BB=\{B_s\}_{s\in S}$ in
$X$, the family $\{\{g\}\times B_s\}_{(g,s)\in G\times S}$ is
uniformly bounded in $G\times X$, so $\{\phi(\{g\}\times
B_s)\}_{(g,s)\in G\times S}$ is uniformly bounded which means
$G\cdot\BB$ is uniformly bounded. \hfill $\blacksquare$

\begin{Rem}
Notice that the infinite dihedral group $Dih_\infty$ acts on
itself by left multiplication so that the action is by uniformly
bornologous functions but the multiplication is not bornologous
(see~\ref{MultBornImpliesBothStructuresTheSame}
and~\ref{MonoFromZtoDih}).
\end{Rem}

\section{Coarsely proper and cobounded actions}\label{section CPactions}

\begin{Def}\label{CoarselyProperDef}
An action $\phi$ of a group $G$ on a coarse space $(X,\CC_X)$ is
{\it coarsely proper} if $\phi_x\colon G\to G\cdot x$ is coarsely
proper for all $x\in X$.
\end{Def}

\begin{Lem} \label{CharOfCoarselyProperActions}
An action $\phi$ of a group $G$ on a coarse space $(X,\CC_X)$ is
coarsely proper if and only if for every bounded subset $U$ of $X$
the family $\{g\cdot U\}_{g\in G}$ is point-finite.
\end{Lem}

\dokaz It follows from the fact $\phi_x^{-1}(U)=\{g\in G\mid x\in
g^{-1}\cdot U\}$ for all $x\in X$ and $U\subset X$. \hfill
$\blacksquare$

\begin{Cor}\label{CPAndUBActions}
If an action $\phi$ of a group $G$ on a coarse space $(X,\CC_X)$
is coarsely proper and by uniformly bornologous functions, then
$\phi_x\colon G\to G\cdot x$ is a coarse equivalence for all $x\in
X$.
\end{Cor}

\dokaz Notice the stabilizer of $x_0$ is finite
by~\ref{CharOfCoarselyProperActions} and use (2)
of~\ref{CreateFirstCS}. \hfill $\blacksquare$

\begin{Lem} \label{CharOfCoarselyProperBornActions}
Let $\phi$ be an action of a group $G$ on a coarse space
$(X,\CC_X)$ by uniformly bornologous functions. Then it is
coarsely proper  if and only if the stabilizer
$$S_U=\{g\in G\mid U\cap(g\cdot U)\ne\emptyset\}$$ of $U$ is
finite for every bounded subset $U$ of $X$.
\end{Lem}

\dokaz One direction is obvious in view of
\ref{CharOfCoarselyProperActions}, so assume $\phi$ is an action
by uniformly bornologous functions that is coarsely proper. If
$S_U=\{g\in G\mid U\cap(g\cdot U)\ne\emptyset\}$ is infinite for
some bounded set $U$, then put $V=\St(U,\{g\cdot U\}_{g\in G}$ and
notice that $\phi^{-1}_x(V)$ contains $S_U$ for all $x\in U$, a
contradiction. \hfill $\blacksquare$

\begin{Def}\label{CoboundedADef}
An action of a group $G$ on a coarse space $(X,\CC_X)$ is {\it
cobounded} if $X=G\cdot U$ for some bounded subset $U$ of $X$.
\end{Def}

\begin{Prop}\label{CharOfCSUnderCoarseA}
If an action $\phi$ of a group $G$ on a coarse space $(X,\CC_X)$
is cobounded and by uniformly bornologous functions, then for
every uniformly bounded family $\BB$ there is a bounded set $U$
such that $\BB$ refines $\{g\cdot U\}_{g\in G}$.
\end{Prop}

\dokaz Pick a bounded set $V$ such that $G\cdot V=X$. Given
$\BB=\{B_s\}_{s\in S}\in \CC_X$ put $U=\St(V,G\cdot \BB)$. $U$ is
bounded and $\BB$ refines $\{g\cdot V\}_{g\in G}$. \hfill
$\blacksquare$

\begin{Cor}\label{CoarseAOnTwoStructures}
If an action $\phi$ of a group $G$ on a set $X$ is cobounded and
by uniformly bornologous functions under two coarse structures
$\CC_1$ and $\CC_2$ on $X$, then $\CC_1=\CC_2$ if and only if
bounded sets in both structures are identical.
\end{Cor}

\dokaz By \ref{CharOfCSUnderCoarseA} both structures are generated
by families $\{g\cdot U\}_{g\in G}$, where $U$ is bounded. \hfill
$\blacksquare$

\section{Coarse actions}\label{section CActions}

\begin{Def}\label{CoarseActionsDef}
An action of a group $G$ on a coarse space $(X,\CC)$ is {\it
coarse} if it is coarsely proper, cobounded, and by uniformly
bornologous functions.
\end{Def}

\begin{Cor}\label{CoarseAInducesCoarseEq}
If an action $\phi$ of a group $G$ on a coarse space $(X,\CC_X)$
is coarse, then $\phi_x\colon (G,\CC_l(G))\to (X,\CC_X)$ is a
coarse equivalence for all $x\in X$.
\end{Cor}

\dokaz By \ref{CPAndUBActions} the function $g\to g\cdot x_0$
is a coarse equivalence from $G$ to $G\cdot x_0$.
Notice the inclusion $G\cdot x_0\to X$ is a coarse equivalence
by the coboundedness of the action.
\hfill $\blacksquare$

\begin{Thm}\label{GenralizedGromovThmTWO}
Suppose $\alpha_i\colon G_i\times X\to X$, $i=1,2$, are two
commutative left actions of groups $G_i$ on the same set $X$. If
there are coarse structures $\CC_i$, $i=1,2$, whose bounded sets
coincide such that $\alpha_i$ is coarse with respect to $\CC_i$,
then
\begin{itemize} \item[a.] $G_1$ is coarsely equivalent to
$G_2$, \item[b.] $(X,\CC_1)$ is coarsely equivalent to
$(X,\CC_2)$.
\end{itemize}
\end{Thm}

\dokaz Pick a bounded set (in both coarse structures) $U$ such
that $G_i\cdot U=X$ for $i=1,2$.  Pick $x_0\in U$. Define
$\psi\colon G_2\to G_1$ so that $h^{-1}\cdot x_0\in \psi(h)\cdot
U$ for all $h\in G_2$.

To show $\psi$ is large scale uniform consider a finite subset $F$
of $G_2$ containing identity, define $V=F ^{-1}\cdot U$ and define
$E$ as the set of all $g\in G_1$ so that $V\cap (g\cdot
V)\ne\emptyset$. Suppose $h=h_1^{-1}h_2\in F$ and $g_i=\psi(h_i)$
for $i=1,2$. Consider $y=g_1^{-1}(h_2^{-1}\cdot x_0)$ and put
$g=g_1^{-1}g_2$. Our goal is to show $y\in V\cap (g\cdot V)$
resulting in $g\in E$. Since $g^{-1}\cdot y=g_2^{-1}(h_2^{-1}\cdot
x_0)\in U\subset V$, $y\in g\cdot V$. Now, as $h_2=h_1\cdot h$,
$y=g_1^{-1}(h_2^{-1}\cdot x_0)= g_1^{-1}(h^{-1}\cdot h_1
^{-1}\cdot x_0)=h^{-1}(g_1^{-1}(h_1 ^{-1}\cdot x_0))\subset
h^{-1}\cdot U\subset F ^{-1}\cdot U=V$.

Similarly, define $\phi\colon G_1\to G_2$ so that $g^{-1}\cdot
x_0\in \phi(g)\cdot U$ for all $g\in G_1$ and notice it is large
scale uniform.

Let $\BB$ be a uniformly bounded family in $\CC_1$ so that all
sets $g\cdot U$, $g\in G_1$, refine $\BB$. Let us observe $g\to
g\cdot x_0$ and $g\to \psi(\phi(g))\cdot x_0$ are
$\St(\BB,\BB)$-close. Indeed, using the definition of $\phi$ and
commutativity of two actions, we get $\phi(g)^{-1}\cdot x_0\in
g\cdot U$, and by definition of $\psi$ we have $\phi(g)^{-1}\cdot
x_0\in \psi(\phi(g))\cdot U$. Since $g\to g\cdot x_0$ is a coarse
equivalence from $G_1$ to $(X_1,\CC_1)$ (see
\ref{CoarseAInducesCoarseEq}), $\psi\circ \phi$ is close to the
identity of $G_1$. Similarly, $\phi\circ \psi$ is close to the
identity of $G_2$. \hfill $\blacksquare$

\section{Topological actions}

Let $X$ be a topological space and $G$ be a group. Recall that an
action of $G$ on $X$ is {\it topologically proper} if each point
$x\in X$ has a neighborhood $U_x$ such that the stabilizer $\{g\in
G\mid U_x\cap (g\cdot U_x)\ne\emptyset\}$ of $U_x$ is finite. An
action of $G$ on $X$ is {\it cocompact} if there exists a compact
subspace $K\subset X$ such that $G\cdot K=X$.

\begin{Def}\label{TopActionsDef}
Let $X$ be a locally compact topological space. An action of a
group $G$ on $X$ is {\it topological} if it is by homeomorphisms,
it is cocompact and topologically proper.
\end{Def}

\begin{Prop}\label{UVW}
Suppose $X$ is a locally compact topological space. If $\phi$ is a
topological action of $G$ on $X$, then there is a unique coarse
structure $\CC_\phi$ on $X$ such that the action $\phi$ of $G$ on
$(X,\CC_\phi)$ is coarse and the bounded sets in $\CC_\phi$ are
precisely relatively compact subsets of $X$. The structure
$\CC_\phi$ is generated by families $\{g\cdot K\}_{g\in G}$ where
$K$ is a compact subset of $X$.
\end{Prop}

\dokaz Uniqueness of $\phi$ follows from
\ref{CoarseAOnTwoStructures}. Let us show that the stabilizer of
each compact subset $K$ of $X$ is finite. If it is not, then there
is an infinite subset $I$ of $G$ and points $x_g\in K\cap (g\cdot
K)$ for each $g\in I$. The set $\{x_g\}_{g\in I}$ must be discrete
(otherwise the action would not be topologically proper at its
accumulation point), so infinitely many $x_g$'s are equal, a
contradiction.

Consider the structure $\CC_\phi$ on $X$ described in the proof
of~\ref{CreateSecondCS}. Notice it has the required properties.
\hfill $\blacksquare$

\begin{Cor}\label{TopActAreEqIfCommutative}
Suppose $\phi\colon G\times X\to X$ and $\psi\colon H\times X\to
X$ are two topological actions on a locally compact space $X$. If
$\phi$ commutes with $\psi$, then $(X,\CC_\phi)$ and
$(X,\CC_\psi)$ are coarsely equivalent.
\end{Cor}

\dokaz Use \ref{GenralizedGromovThmTWO}.
\hfill $\blacksquare$

\begin{Rem}\label{TopActAreMayBeDifferentRem}
It is not true that $\CC_\phi=\CC_\psi$ in general. Use
\ref{TwoCSExistOnDih} and equip groups with discrete topologies.
\end{Rem}

\begin{Thm}\label{GenralizedGromovThm}
If $G$ and $H$ are coarsely equivalent groups, then there is a
locally compact topological space $X$ and topological actions
$\phi\colon G\times X\to X$ and $\psi\colon H\times X\to X$ that
commute.
\end{Thm}

\dokaz Pick a coarse equivalence $\alpha\colon G\to H$. Choose a
function $c$ assigning to each finite subset $F$ of $G$ a finite
subset $c(F)$ of $H$ with the property that $u^{-1}\cdot v\in F$
implies $\alpha(u)^{-1}\cdot \alpha(v)\in c(F)$.

Choose a function $d$ assigning to each finite subset $F$ of $H$
a finite subset $d(F)$ of $G$ with the property that
$\alpha(u)^{-1}\cdot \alpha(v)\in F$ implies $u^{-1}\cdot v\in
d(F)$.

Let $E$ be a finite subset of $H$ so that $H=\alpha(G)\cdot E$.

Let $X$ be the space of all functions $\beta\colon G\to H$
satisfying the following conditions:
\begin{enumerate}
\item $u^{-1}\cdot v\in F$
implies $\beta(u)^{-1}\cdot \beta(v)\in c(F)$ for all finite subsets $F$ of $G$,
\item $\beta(u)^{-1}\cdot \beta(v)\in F$
implies $u^{-1}\cdot v\in d(F)$ for all finite subsets $F$ of $H$.
\item $H=\beta(G)\cdot E$.
\end{enumerate}
We consider $X$ with the compact-open topology provided both $G$
and $H$ are given the discrete topologies. Notice $X$ is closed in
the space $H^G$ of all functions from $G$ to $H$ equipped with the
compact-open topology. Indeed, Conditions (1) and (2) above hold
for all $\beta\in cl(X)$, so it remains to check $H=\beta(G)\cdot
E$ for such $\beta$. Given $h\in H$ consider the set
$F=\beta(1_G)^{-1}\cdot h\cdot E^{-1}$ and choose $\gamma\in X$ so
that $\gamma(g)=\beta(g)$ for all $g\in d(F)\cup\{1_G\}$. Pick
$g_1\in G$ and $e\in E$ so that $h=\gamma(g_1)\cdot e$. Since
$\gamma(1_G)^{-1}\cdot \gamma (g_1)\in F$, $g_1=1_G^{-1}\cdot
g_1\in d(F)$ and $\gamma(g_1)=\beta(g_1)$. Thus $h\in
\beta(G)\cdot E$.

Notice $X$ is locally compact. Indeed, given $\beta\in X$ consider
$U=\{\gamma\in X\mid \gamma(1_G)=\beta(1_G)\}$. It is clearly open
and equals $X\cap K$, where $K\subset H^G$ is the set of all
functions $u$ satisfying $u(g)\in \beta(1_G)\cdot c(\{g\})$.
Notice $K$ is compact (it is a product of finite sets). Since $X$
is closed in $H^G$, $X\cap K$ is compact as well.

The action of $G$ on $X$ is given by $(g\cdot \beta)(x):=\beta (g\cdot x)$.
The action of $H$ on $X$ is given by $(h\cdot \beta)(x):=h\cdot \beta (x)$.
Notice that the two actions commute.
The action of $H$ on $X$ is cocompact: $X=H\cdot K$,
where $K=\{\beta\in X\mid \beta(1_G)=1_H\}$.
The action of $G$ on $X$ is cocompact: $X=G\cdot L$,
where $L$ is the set of $\beta\in X$ such that $\beta(1_G)\in E^{-1}$
(which implies $\beta(g)\in E^{-1}\cdot c(\{g\})$ for all $g\in G$ so that $L$
is compact). Indeed, for any $\gamma\in X$ there is $e\in E$
such that $1_H=\gamma(g_1)\cdot e$ for some $g_1\in G$.
Put $\beta(x)=\gamma(g_1\cdot x)$ and notice $\beta(1_G)=e^{-1}\in E^{-1}$,
so $\beta\in L$ and $\gamma=g_1\cdot\beta$.

Action of $H$ is proper: given $\beta\in X$ put $U=\{\gamma\in
X\mid \gamma(1_G)=\beta(1_G)\}$. If $\lambda\in U\cap (h\cdot U)$,
then $\lambda(1_G)=\beta(1_G)$ and $h^{-1}\cdot
\lambda(1_G)=\beta(1_G)$. Thus $h=1_H$.

Action of $G$ is proper: given $\beta\in X$ put $U=\{\gamma\in
X\mid \gamma(1_G)=\beta(1_G)\}$. If $\lambda\in U\cap (g\cdot U)$,
then $\lambda(1_G)=\beta(1_G)$ and $\lambda(g^{-1})=\beta(1_G)$.
Thus $\lambda(g^{-1})=\lambda(1_G)$ which implies $g^{-1}\in
d(\{1_H\})$, so the set of such $g$ is finite. \hfill
$\blacksquare$

\end{document}